\documentclass[11pt]{article}
\usepackage{amssymb,amsmath,latexsym,bbm,amscd}

\newtheorem{example}[equation]{Example}
\newtheorem{theorem}[equation]{Theorem} 
\newtheorem{corollary}[equation]{Corollary}
\newtheorem{lemma}[equation]{Lemma}
\newtheorem{proposition}[equation]{Proposition}
 
\newtheorem{remark}[equation]{Remark} 
\numberwithin{equation}{section}

\newcommand{\Cx}{{\mathbb C}}

\newcommand{\g}{\mathfrak{g}}

\newcommand{\ot}{\otimes}
\newcommand{\ol}{\overline}
\newcommand{\Z}{\mathbb{Z}}

\newcommand{\gs}{\sigma}

\newcommand{\proof}{{\bf Proof\ \ }}
\newcommand{\qed}{\hfill $\Box$}

\newcommand{\A}{{\mathcal A}}
\newcommand{\B}{{\mathcal B}}
\newcommand{\End}{\hbox{End\,}}
\newcommand{\M}{\widetilde{M}}

\newcommand{\bm}{{\bf m}}
\newcommand{\bt}{{\bf t}}
\newcommand{\C}{{\mathbb C}}
\newcommand{\K}{{\mathbb C}}
\newcommand{\R}{{\mathcal R}}
\newcommand{\CC}{{\mathcal C}}
\newcommand{\ZZ}{{\mathcal Z}}

\newcommand{\Aut}{{\rm Aut}}
\newcommand{\GL}{{\rm GL}}
\newcommand{\Hom}{{\rm Hom}}
\newcommand{\id}{{\rm id}}
\newcommand{\Id}{{\rm id}}
\newcommand{\Span}{{\rm Span}}

\newcommand{\wG}{{\widehat G}}
\newcommand{\wH}{{\widehat H}}

\newcommand{\hro}{{\widehat \rho}}
\newcommand{\hM}{{\widehat M}}
\newcommand{\tN}{{\widetilde N}}
\newcommand{\ba}{{\bf a}}
\newcommand{\bb}{{\bf b}}
\newcommand{\T}{{\bf T}}

\title{Thin Coverings of Modules}

\date{}

\author{ Yuly Billig$\hbox{\,}^1\hbox{}^*$\ \, and Michael Lau$\hbox{\,}^2$\thanks{Both authors gratefully acknowledge funding 
from the Natural Sciences and Engineering Research Council of Canada.  
Y.B. is partially supported by an NSERC Discovery Grant; M.L. is supported by an NSERC postdoctoral fellowship.}\ \,\footnote{Corresponding author.} \vspace{0.3cm}\\$\hbox{\ \,}^1${\small Carleton University, School of Mathematics and Statistics},\\ {\small 1125 Colonel By Drive, Ottawa, Ontario, Canada K1S 5B6}\\ {\small Email:\ billig@math.carleton.ca}\vspace{0.1cm}\\
$\hbox{\ \,}^2${\small University of Ottawa, Department of Mathematics and Statistics,}\\{\small 585 King Edward Ave., Ottawa, Ontario, Canada K1N 6N5}\\{\small  Email:\ mlau@uottawa.ca}}

\begin{document}
\maketitle

\begin{small}
\noindent 
{\bf Abstract.} Thin coverings are a method of constructing graded-simple modules from simple (ungraded) modules.  
After a general discussion, we classify the thin coverings of (quasifinite) simple modules over 
associative algebras graded by finite abelian groups.  The classification uses the representation theory of 
cyclotomic quantum tori. We close with an application to representations of multiloop Lie algebras.

\bigskip



%\noindent
%{\bf Keywords: }

\end{small}

\setcounter{section}{-1}
\vskip.25truein
\section{Introduction}

In a recent series of papers \cite{ABP1,ABP2,ABFP}, B.~Allison, S.~Berman, A.~Pianzola, and J.~Faulkner 
examined 
the structure of multiloop algebras.  These algebras are formed by a generalization of the twisting process 
used in Kac's construction \cite{kac} of the affine Lie algebras.

More specifically, for any algebra $\B$ over the complex field $\Cx$ and finite-order commuting 
automorphisms $\gs_1,\ldots,\gs_N$ of $\B$ of period $m_1,\ldots, m_N$, respectively, the multiloop 
algebra $L(\B,\gs_1,\ldots,\gs_N)$ is the following subalgebra of $\B\ot\Cx[t_1^\pm,\ldots,t_N^\pm]$:
$$L(\B,\gs_1,\ldots,\gs_N)=\bigoplus_{(i_1,\ldots,i_N)\in\Z^{N}}\B_{\ol{\iota}_1,\ldots,\ol{\iota}_N}\ot
\Cx t_1^{i_1}\cdots t_N^{i_N},$$
where $\B_{\ol{\iota}_1,\ldots,\ol{\iota}_N}$ is the common eigenspace
$$\{b\in\B\ |\ \gs_j b =\xi_j^{i_j} b \hbox{\ for \ }1\leq j\leq N\}$$
for fixed primitive $m_j$th roots of unity $\xi_j$.

In particular, when $\B$ is a finite-dimensional simple Lie algebra, one can construct realizations for almost 
all extended affine Lie algebras by taking central extensions of the resulting multiloop algebras and adjoining 
appropriate derivations \cite{ABFP, neher}.

While studying the representation theory of such extended affine Lie algebras 
\cite{EALAreps}, we considered methods of ``twisting'' simple $\B$-modules 
into graded-simple modules for the multiloop algebras  $L(\B,\gs_1,\ldots,\gs_N)$.  
As part of our construction, we used structures called thin coverings of modules.  

Let $\B=\bigoplus_{g\in G}\B_g$ be an algebra graded by an abelian group $G$.  A {\em thin covering} of a left 
$\B$-module $M$ is a family of subspaces $\{M_g\ |\ g\in G\}$ with the following properties:

(i)\quad\ \ $M=\sum_{g\in G}M_g$,

(ii)\quad\ $\B_gM_h\subseteq M_{g+h}$ for all $g,h\in G$,

(iii)\quad If $\{ N_g\ |\ g\in G\}$ also satisfies (i) and (ii) and $N_g\subseteq M_g$ for all $g\in G$, 
then $N_g=M_g$ for all $g\in G$.

We consider finite-dimensional and infinite-dimensional quasifinite modules over a $G$-graded 
associative algebra $\A$.  
%While defined axiomatically, these structures appear in graded Clifford theory for associative algebras, 
%relating the categories of left $\A_0$-modules and graded left $\A$-modules (cf.~\cite{NaOy}, for instance).  
Each thin covering of an (ungraded) left $\A$-module $M$ is associated with a graded left $\A$-module $\M$.  
Although thin coverings are generally not unique, the isomorphism class of $M$ completely determines 
the graded-isomorphism class of $\M$, under mild conditions.  It is trivial to show that $\M$ is graded-simple 
whenever $M$ is simple.

An $\A$-module may be twisted with an $\A$-automorphism associated with its $G$-grading.  Such twists do not 
change the isomorphism class of the graded modules coming from the thin covering construction, and they play 
a vital role in our classification of thin coverings.

After these applications to graded theory, we give an explicit characterization of thin coverings of 
quasifinite simple modules over arbitrary (unital) associative algebras graded by finite abelian groups.  
The classification uses the fact that any such module is also a module for a cyclotomic quantum torus coming 
from isomorphisms between twists of the module.

We conclude the paper with an application to the representation theory of multiloop Lie algebras.  Our approach is an alternative to Clifford theory, where graded-simple $\A$-modules are constructed by 
inducing from simple $\A_0$-modules.  (See \cite{NaOy}, Thm.~2.7.2 and also \cite{clifford, dade}.)
 
This method gives an explicit procedure for constructing graded-simple $\A$-modules from (ungraded) simple 
$\A$-modules via the action of a cyclotomic quantum torus.
For the applications we have in mind, it would be difficult to get a classification of simple $\A_0$-modules,
whereas ungraded simple $\A$-modules are well-understood.

\medskip

\noindent
{\bf Note:} Although we work over the field $\C$ of complex numbers, all the results of this paper hold over 
an arbitrary algebraically closed field of characteristic zero.  All modules are left modules unless otherwise indicated.  Likewise, all hypotheses of simplicity and semisimplicity should be interpreted as left simplicity and left semisimplicity.

\medskip

\noindent
{\bf Acknowledgments.}
Y.B.~would like to thank the University of Sydney, and M.L.~the University of Alberta, for their warm 
hospitality during our visits.  Y.B.~also expresses grateful appreciation for helpful conversations with 
C.~Dong.

\section{Thin Coverings and Graded Modules}\label{one}

%In this section, we introduce {\em thin coverings}, a way of associating graded-simple modules with simple 
%(not necessarily graded) modules over graded associative algebras.  We classify thin coverings over irreducible 
%finite-dimensional modules; then we apply the general finite-dimensional theory to the locally finite case of 
%highest weight modules for affine Lie algebras.  The process will later be used to twist modules for the 
%toroidal Lie algebras into modules for the corresponding twisted toroidal Lie algebras.

Let $\A$ be a (unital) associative algebra graded by an abelian group $G$.  For any
module $M$ over $\A$, a set of subspaces $\{M_g\ |\ g\in G\}$ is a $G$-{\em covering} of $M$ 
%(relative to $(\A, G)$) 
if $\sum_{g\in G}M_g=M$ and $\A_gM_h\subseteq M_{g+h}$ for all $g,h\in G$.  
We will drop the prefix $G$ from ``$G$-covering'' when there is no ambiguity.
Two coverings $\{M_g\ |\ g\in G\}$ and $\{M_g'\ |\ g\in G\}$ are {\em equivalent}, and we write
$$\{M_g\ |\ g\in G\}\sim \{M_g'\ |\ g\in G\}$$
if there is a fixed $h\in G$ such that $M_g=M_{g+h}'$ for all $g\in G$.  

The set $\mathfrak{C}:={\mathfrak C}(M,\A,G)$ of coverings of the module $M$ is 
partially ordered:  for any $\{M_g\ |\ g\in G\}$ and $\{N_g\ |\ g\in G\}$ in $\mathfrak{C}$, let
$$\{M_g\ |\ g\in G\}\preceq \{N_g\ |\ g\in G\}$$
if the covering $\{M_g\}$ is equivalent to a covering $\{M_g'\}$ for which $M_g'\subseteq N_g$ for all $g\in G$.

The minimal coverings in the poset $\mathfrak{C}$ are called {\em thin coverings}.  
 If $M$ is a finite-dimensional or infinite-dimensional quasifinite module (see Section 4), then it has
a thin covering.

\begin{remark}
{\em Let $\g$ be a Lie algebra graded by an abelian group $G$, and let $M$ be a $\g$-module. The grading on $\g$
extends naturally to a $G$-grading on the universal enveloping algebra $U(\g)$. Note that the coverings of $M$
as a $\g$-module coincide with its coverings as a $U(\g)$-module. Thus, by studying thin coverings of modules
over associative algebras, we automatically obtain the corresponding results for Lie algebras.
The same applies to Lie superalgebras as well, where we would need to consider the $\Z_2$-graded version of 
a thin covering.  (See the definition in Section 4.)}
\end{remark}

 The following elementary, but important lemma says that the thin coverings 
of a simple $\A$-module $M$ are essentially determined by the simple $\A_0$-submodules of $M$.

\medskip

\begin{lemma}\label{obvious}
Let $\{M_g\ |\ g\in G\}$ be a covering for a simple module $M$ over an associative algebra $\A$ graded by an 
abelian group $G$.  
This covering is thin if and only if 
$$M_g=\A_{g-h}u,$$
for all $g,h\in G$ and any nonzero $u\in M_h$.

In particular, in a thin covering, every nonzero space $M_h$ is a simple $\A_0$-module.
\end{lemma}
\proof
Since $M$ is a simple $\A$-module, $\sum_{g\in G}\A_{g-h}u=\A u=M$.  
This gives a covering of $M$ with $\A_{g-h}u\subseteq M_g$ for all $g\in G$.  
Thus $\{ M_g\}$ is a thin covering if and only if $\A_{g-h}u=M_g$ for all $g$.\qed

\bigskip

Although the sum $\sum_{g\in G}M_g$ of the spaces $M_g$ in a thin covering need not be direct, the {\em external}
direct sum of these subspaces always has a natural graded module structure.  In fact, this construction 
``preserves simplicity'', as described in the following theorem. 

\begin{theorem}\label{grmodlm} Let $M$ be a simple (not necessarily graded) module over an associative 
algebra $\A$ graded by an abelian group $G$.  Then for any thin covering $\{M_g\ |\ g\in G\}$ of $M$, the space 
$$\widetilde{M}:=\bigoplus_{g\in G}M_g$$
is a graded-simple $\A$-module.
\end{theorem}
\proof
The space $\widetilde{M}$ has an obvious graded module structure, namely $\widetilde{M}_g:=M_g$ and 
$\A_g:\ \widetilde{M}_h\rightarrow\widetilde{M}_{g+h}$ for all $g,h\in G$.

By Lemma \ref{obvious}, $\A_{h-g}u=M_h$ for all $g\in G$ and nonzero elements $u\in M_g$.  Likewise, 
$A_{h-g}\widetilde{u}=\widetilde{M}_h$ for every $h\in G$, where $\widetilde{u}$ is the element of $\M$ having 
the nonzero element $u\in M_g$ in the $\M_g$-component and $0$ elsewhere.

Suppose $\widetilde{N}=\bigoplus_{g\in G}\widetilde{N}_g$ is a nonzero $G$-graded submodule of $\M$.  
Then for any nonzero component $\widetilde{N}_g$ and $h\in G$, we have $\A_{h-g}\widetilde{N}_g=\M_h$.  
Hence $\widetilde{N}=\M$.\qed

\medskip

In the next theorem, we show that graded-simple $\A$-modules come from thin coverings 
of simple (non-graded) $\A$-modules. The assumption on the existence of a maximal submodule holds
in the finite-dimensional case as well as for the infinite-dimensional quasifinite modules.  (See Remark
\ref{qmax}.)  

\begin{theorem} \label{converse} 
Let $\M$ be a graded-simple module over an associative algebra
$\A$ graded by an abelian group $G$:
$$ \M = \mathop\bigoplus\limits_{g\in G} \M_g .$$
Assume that as a non-graded $\A$-module, $\M$ has a maximal (non-graded) 
submodule and a simple (non-graded) quotient $M$:
$$ \varphi:  \ \M \rightarrow M .$$
Let $M_g = \varphi(\M_g)$.
Then $\{ M_g\ |\ g\in G \}$ is a thin covering of the module $M$ and
$\M$ is graded-isomorphic to the graded $\A$-module associated with this thin covering,
$$ \M \cong \mathop\bigoplus\limits_{g\in G} M_g ,$$
as in Theorem \ref{grmodlm}.
\end{theorem}
\proof
Since the map $\varphi$ is surjective, we get that 
$\sum\limits_{g\in G} \varphi(\M_g) = \varphi(\M) = M $,
and hence $\{ \varphi(\M_g)\ |\ g\in G \}$ is a covering of $M$.

Moreover, the restrictions of $\varphi$ to $\M_g$,
$$\varphi_g: \M_g \rightarrow M_g ,$$
are bijections. If these maps had non-trivial kernels, then
$\mathop\bigoplus\limits_{g\in G} \ker\, \varphi_g$ would be a non-trivial
graded submodule of $\M$, which contradicts the fact that $\M$ is graded-simple. 

We now show that the covering $\{ \varphi(\M_g)\ |\ g\in G \}$ is thin. 
Suppose $\{ N_g\ |\ g\in G \}$ is another covering of $M$ with $N_g \subseteq M_g$.
Let $\tN_g \subseteq \M_g$ be the pre-image of $N_g$ under the map $\varphi_g$. 
It is easy to see that
$$\tN = \mathop\bigoplus\limits_{g\in G} \tN_g $$
is a graded $\A$-submodule in $\M$. To verify this, it is enough to check that
$\A_h \tN_g \subseteq \tN_{g+h}$. However, $\varphi (\A_h \tN_g) =
\A_h N_g \subseteq N_{g+h}$, and $\tN_{g+h}$ is the pre-image of $N_{g+h}$ under
$\varphi_{g+h}$. Thus $\A_h \tN_g \subseteq \tN_{g+h}$, and $\tN$ is a nonzero graded
submodule of $\M$. Since $\M$ is graded-simple, we conclude that
$\tN_g = \M_g$ for all $g\in G$. This implies that $N_g = M_g$ because the map
$\varphi_g$ is a bijection. Hence the covering $\{ M_g\ |\ g\in G \}$ is indeed thin.

Once again using bijectivity of $\varphi_g$, we see that the graded $\A$-module
associated with this covering is isomorphic to $\M$:
$$ \mathop\bigoplus\limits_{g\in G} M_g \cong  \mathop\bigoplus\limits_{g\in G} \M_g.$$
\qed

\bigskip

Gradings on $\A$ by finite abelian groups can be alternatively described via finite
abelian subgroups of the group of automorphisms of $\A$. 
Suppose $\A$ is graded by a finite abelian group $G$, and consider the dual group $\wG = \Hom (G, \C^*)$.
Then we can interpret an element $\sigma\in \wG$ as the finite order $\A$-automorphism
defined by $\sigma(a) = \sigma(g) a$ for $a\in \A_g$.
Although the groups $G$ and $\wG$ are isomorphic, there is no canonical isomorphism
between them.

Conversely, a finite abelian group $\wG$ of automorphisms of $\A$ defines a grading. 
The algebra $\A$ is graded by the group $G = \Hom (\wG, \C^*)$:
$$\A=\bigoplus_{g\in G}\A_{g},$$
where $\A_g = \left\{ a\in \A\ \left|\ \sigma(a) = g(\sigma) a, \text{\rm \ for all \ }
\sigma\in\wG \right.\right\}$. 

To be consistent with the additive
notation for grading groups, we will slightly abuse notation and always treat $G$ as an additive
group.  The sum of two elements $g,h\in G$ is defined by 
\begin{equation}\label{addnotation}
(g+h)(\sigma) = 
g(\sigma) h(\sigma)\hbox{\ for\ }\sigma\in\wG,
\end{equation}
 and we will denote the identity element of $G$ by 0.

An $\A$-module $(M,\rho)$ can be twisted by an automorphism $\gs \in \Aut (\A)$, and we use the notation
$M^{\gs}$ for the module $(M,\rho\gs)$.  That is, the module action in $M^{\gs}$ 
is given by $\rho\gs:\ \A\rightarrow \End_\Cx M$.

The following proposition shows that a simple $\A$-module and its twists with automorphisms in $\wG$ 
yield isomorphic graded-simple modules.

\begin{proposition}\label{twists}
Let $\A$ be an associative algebra graded by a finite abelian group $\wG$.  Let $M$ be a simple 
 $\A$-module, and let $\sigma\in\wG$.

(i)\ \ Any thin covering $\{ M_g\ |\ g\in G \}$ of $M$ is also a thin covering of $M^\sigma$.

(ii)\ Let $\{ M_g\ |\ g\in G \}$ be a thin covering of $M$, and let $\{ M_g^\sigma\ |\ g\in G \}$
be the same covering of $M^\sigma$. Then the following graded-simple  $\A$-modules are naturally isomorphic
$$ \mathop\bigoplus\limits_{g\in G} M_g^\sigma \cong \left( \mathop\bigoplus\limits_{g\in G} 
M_g \right)^\sigma .$$

(iii) If $M$ is a $G$-graded  $\A$-module, then $M^\sigma \cong M$ for all $\gs\in \wG$.
\end{proposition}
\proof
We first show that a covering of $M$ is also a covering of $M^\sigma$. For all $a\in A_g$ and $m\in M_h$, we 
have
$\rho(a) m \in M_{g+h}$, so 
$$\rho \sigma(a) m = g(\sigma) \rho(a) m \in M_{g+h}.$$
This gives a bijection between the coverings of $M$ and the coverings of $M^\sigma$. Clearly the thin coverings 
of $M$
correspond to the thin coverings of $M^\sigma$.

Part (ii) is obvious. 
 
 To establish (iii), we define the isomorphism 
$$\theta: M = \mathop\bigoplus\limits_{g\in G} M_g \rightarrow M^\sigma$$
by $\theta(m) = g(\sigma) m$ for $m\in M_g$. Let us check that $\theta$ commutes with the action of $\A$.
For $a \in \A_g, m\in M_h$, we have
$$\rho \sigma(a) \theta (m) = g(\sigma) h(\sigma) \rho(a) m = \theta( \rho(a) m) .$$ 
\qed
 
\begin{remark}
 {\em The following partial converse is true: if the modules $M$ and $M^\sigma$
are isomorphic for some finite order automorphism $\gs\in\Aut\A$, then $M$ admits a grading by the cyclic group 
$\left< \sigma \right>$.  (See Example \ref{cycex} below.) However, if we replace a cyclic group 
$\left< \sigma \right>$ by a finite
abelian group $G$, the analogous statement is, in general, false (cf.~Example \ref{matrixex}).}  
\end{remark}

\section{Thin coverings of modules over a semisimple algebra}

Recall that an associative algebra $\A$ is {\em  semisimple} if each  module for $\A$ is completely 
reducible.  In Theorem \ref{unicite}, we prove that if $\A$ is semisimple, $1\in\A_g\A_{-g}$ for all $g$, and 
each $\A_g$ is an irreducible  $\A_0$-module, then the isomorphism class of the graded-simple module 
$\widetilde{M}$ does not depend on the choice of the thin covering of $M$.  The following lemma is well-known 
(cf.~\cite{NaOy}, Cor. 2.74, for instance):

\begin{lemma}\label{sslm}
Let $\A$ be a semisimple associative algebra.  If $\A$ is graded by an abelian group $G$, then its 
subalgebra $\A_0$ is also semisimple.\qed
\end{lemma}

\begin{theorem}\label{unicite}
Let $\A$ be a semisimple associative algebra graded by an abelian group $G$.  Assume that for each 
$g\in G$, the graded component $\A_g$ satisfies the following two properties:

(i)\quad $1\in\A_g\A_{-g}$,

(ii)\quad $\A_g$ is irreducible as a  $\A_0$-module.

\noindent
Let $\{M_g\ |\ g\in G\}$ and $\{N_g\ |\ g\in G\}$ be thin coverings for an irreducible  $\A$-module $M$.  
Then for some $h\in G$,
$$\bigoplus_{g\in G}M_g\cong\bigoplus_{g\in G}N_{g-h}.$$  
as graded  $\A$-modules.
\end{theorem} 
\proof If $M=0$, there is nothing to prove.  Otherwise, choose $h\in G$ so that $N_h\neq 0$.  Note that every 
$N_g$ is nonzero, since 
$$0\neq N_h\subseteq \A_{h-g}\A_{g-h}N_h\subseteq \A_{h-g}N_g .$$ 
% by Lemma \ref{obvious}.  
Likewise, $M_g\neq 0$ for all $g\in G$.  By Lemma \ref{sslm}, $N_0$ has a direct complement $H$ as an 
$\A_0$-module.  Since $\{M_g\ |\ g\in G\}$ is a covering, there is an $h\in G$ such that the projection
$$\pi_h:\ M_h\rightarrow N_0$$
(given by the splitting $M_h\subseteq M=N_0\oplus H$) is nonzero.  Since $\pi_h$ is a nonzero homomorphism of 
simple $\A_0$-modules (using Lemma \ref{obvious}), $M_h\cong_{\A_0}N_0$ by Schur's Lemma.  

Our next task is to show that 
$$\bigoplus_{g\in G}M_g\cong\A\otimes_{\A_0}M_h$$
as graded  $\A$-modules.  We begin by noting that for each $g\in G$, the  $\A_0$-module 
$\A_g\otimes_{\A_0}M_h$ is irreducible.

Let $0\neq\sum_ia_i\ot m_i\in\A_g\ot_{\A_0}M_h$, and let $m\in M_h$ be nonzero.  Since $M_h$ is irreducible as 
an $\A_0$-module (Lemma \ref{obvious}), there exist $b_i\in\A_0$ such that $b_i m=m_i$ for all $i$.  Then 
$$\sum_ia_i\ot m_i=\left(\sum_ia_ib_i\right)\ot m.$$
Since $\A_g$ is irreducible as a  $\A_0$-module, we have
\begin{eqnarray*}
\A_0 \left(\left(\sum_ia_ib_i\right)\ot m\right)&=&\A_g\ot m\\
&\subseteq&\A_g\A_0\ot m\\
&=&\A_g\ot\A_0m\\
&=&\A_g\ot_{\A_0} M_h.
\end{eqnarray*}
Thus $\A_g\ot_{\A_0}M_h$ (and $\A_g\ot_{\A_0}N_0$) is an irreducible  $\A_0$-module.

Note that $\A\ot_{\A_0}M_h=\bigoplus_{g\in G}\left(\A_g\ot_{\A_0}M_h\right)$, and there is a well-defined graded 
map
$$\phi:=\bigoplus_{g\in G}\phi_g:\ \bigoplus_{g\in G}\left(\A_g\ot_{\A_0}M_h\right)\rightarrow 
\bigoplus_{g\in G}M_g$$
given by $\phi_g(a\ot m)=am\in M_g$ for all $a\in \A_{g-h}$ and $m\in M_h.$  The map $\phi_g$ is a nonzero 
homomorphism between simple $\A_0$-modules, so is an $\A_0$-module isomorphism.  Thus $\phi$ is a bijection.  
Also 
\begin{eqnarray*}
\phi(ab\ot m)&=&abm\\
&=&a\phi(b\ot m)
\end{eqnarray*}
for any homogeneous $a,b\in\A$ and $m\in M_h$, so $\phi$ is an isomorphism of graded  $\A$-modules.  
Likewise,
$$\A\ot_{\A_0}N_0=\bigoplus_{g\in G}\A_g\ot_{\A_0}N_0\cong\bigoplus_{g\in G}N_g$$
as graded  $\A$-modules.

Since the $\A_0$-modules $M_h$ and $N_0$ are isomorphic, we get that $\A \otimes_{\A_0} M_h \cong \A \otimes_{\A_0}
N_0$, and thus
$$\bigoplus_{g\in G}M_g\cong\bigoplus_{g\in G}N_{g-h}.$$  
as graded  $\A$-modules.
\qed

\bigskip
\bigskip

While this method of associating graded modules $\widetilde{M}$ with modules $M$ is not functorial (it is not 
even well-defined if $M$ is not simple), it is surjective in the sense that not only does every graded module 
come from a thin covering, but (in the context of the previous theorem) every graded-simple module comes from a 
thin covering of a simple module.  

The first assertion (``graded modules come from thin coverings'') is trivial.  Let $M=\bigoplus_{g\in G}M_g$ 
be a graded module of a graded associative algebra $\A=\bigoplus_{g\in G}\A_g$.  Forgetting the graded 
structure on $M$ gives an (ungraded) $\A$-module $M$ with thin covering $\{ M_g\ |\ g\in G\}$.  The second 
claim (``graded-simples come from simples'') is the following theorem.  (See also Theorem \ref{converse}.)

\begin{theorem} \label{completeness} Let $M=\bigoplus_{g\in G}M_g$ be a graded-simple module for a 
semisimple associative algebra $\A=\bigoplus_{g\in G}\A_g$ graded by an abelian group $G$.  Assume that 
$\A$ satisfies conditions (i) and (ii) of Theorem \ref{unicite}.  Then there is a simple module $N$ 
for $\A$ so that $N$ has a thin covering $\{ N_g\ |\ g\in G\}$ with
$$\bigoplus_{g\in G}N_g\cong\bigoplus_{g\in G}M_g$$
as graded  $\A$-modules.
\end{theorem}
\proof Let $N$ be a simple  $\A$-submodule of $M$, and let $N'$ be a simple  $\A_0$-submodule of $N$.  
Since $N$ is simple, 
$$N=\A N'= \sum_{g\in G}\A_gN'.$$
Clearly $\A_g\A_hN'\subseteq\A_{g+h}N'$, so $\{\A_gN'\ |\ g\in G\}$ is a covering of $N$.  By the argument in 
the proof of Theorem \ref{unicite}, the $\A_0$-modules $\A_gN'$ are simple, as are the components $N_g$ in 
any subcovering $\{N_g\ |\ g\in G\}$.  The $N_g$ are all nonzero by the proof in Theorem \ref{unicite}.  
Thus the covering $\{\A_gN'\ |\ g\in G\}$ is thin.

Since $N$ is nonzero, the projection 
$$\pi_h:\ N\rightarrow M_h$$ 
is nonzero for some $h\in G$.  Choose $h'\in G$ so that the restriction of $\pi_h$ to $\A_{h'}N'$ is also 
nonzero.  For each $g\in G$, let 
$$C_g=\A_{g-h+h'}N'.$$
Then $\{C_g\ |\ g\in G\}$ is a thin covering of $N$ with the property that the projection map
$$\pi_h:\ C_h\rightarrow M_h$$
is nonzero.  Since $C_g=\A_{g-h}C_h$ and $M_g=\A_{g-h}M_h$, we see that the projection 
$\pi_g:\ C_g\rightarrow M_g$ is also nonzero.  

Let 
$$\pi=\bigoplus_{g\in G}\pi_g:\ \bigoplus_{g\in G}C_g\rightarrow \bigoplus_{g\in G}M_g,$$
where the restriction of $\pi$ to $C_g$ is the projection $\pi_g$.  
Then since each $\pi_g$ is a nonzero homomorphism between simple $\A_0$-modules $C_g$ and $M_g$, the maps 
$\pi_g$ 
are $\A_0$-module isomorphisms.  Thus $\pi$ is a bijection.

Finally, $\pi$ is also a graded homomorphism since
$$\pi(an)=\pi_{g+g'}(an)=a\pi_{g'}(n)=a\pi(n)$$
for all $a\in \A_g$ and $n\in C_{g'}$.  Thus 
$$\bigoplus_{g\in G}C_g\cong\bigoplus_{g\in G}M_g$$
as graded  $\A$-modules.\qed

\section{Classification of Thin Coverings}

%\bigskip

In this section, we will assume that $\A$ is a (unital) associative algebra over $\Cx$, and $M$ is a 
finite-dimensional irreducible $\A$-module whose action is given by the homomorphism 
$\rho:\ \A\rightarrow \End_\Cx M.$

 In order to find the thin coverings of a module $M$, it is important to know for which
$\gs\in\wG$, the modules $M$ and $M^\gs$ are isomorphic.

\medskip

\begin{lemma}\label{Hsubgroup}
Let $\wG$ be a subgroup of $\Aut\,\A$ and let
$$\wH = \left\{\left. \gs\in\wG\ \right|\ M^\gs \cong M \text{\rm \ as
\ } \A \text{\rm -modules} \right\}.$$
Then

(i)\ \  $\wH$ is a subgroup of $\wG$.

(ii)\  For $\gs_1, \gs_2 \in\wG$ we have $M^{\gs_1} \cong M^{\gs_2}$ if and only if
$\gs_1 \gs_2^{-1} \in \wH$.
\end{lemma}
\proof
An isomorphism between $\A$-modules $M$ and $M^\sigma$ is a map $T_\gs \in \GL (M)$
satisfying
$$ T_\gs \rho(a) u = \rho(\sigma(a)) T_\gs u ,$$
for all $a\in\A$ and $u\in M$. Equivalently,
$$ \rho(\sigma(a)) = T_\gs \rho(a) T_\gs^{-1} .$$
The first part of the lemma follows from this relation.

The fact that the modules $M^{\gs_1}$ and $M^{\gs_2}$ are isomorphic is equivalent
to the existence of a $T\in \GL (M)$, such that 
$\rho(\sigma_1 (a)) = T \rho(\sigma_2 (a)) T^{-1}$. Letting $b = \gs_2 (a)$, we get 
that $\rho(\gs_1 \gs_2^{-1} (b)) = T \rho(b) T^{-1}$, which means that 
$\gs_1 \gs_2^{-1} \in\wH$. This establishes the second claim of the lemma.
\qed

\bigskip
\bigskip

 Consider the operator $T_\gs$ introduced in the above proof. If $\gs\in\wH$ has order
$m$, then for all $a\in\A$,
$$\rho(a) = \rho(\gs^m (a)) = T_\gs^m \rho(a) T_\gs^{-m} .$$
By Schur's Lemma, $T_\gs^m$ is a scalar operator, and we may normalize $T_\gs$ to get
$T_\gs^m = \id$. 
It is often convenient to use this normalization of $T_\gs$.

Recall that there is an equivalent description of gradings and coverings in terms of the dual groups 
$H=\Hom(\wH,\Cx^*)$ and $G=\Hom(\wG,\Cx^*)$ (cf.~Section \ref{one}).  Our strategy will be to describe the thin 
$H$-coverings of $M$, and then use them to 
determine the thin $G$-coverings. In order to understand the structure of thin $H$-coverings,
we need to better understand the $H$-grading on $\A$ and on its image $\rho(\A) = \End_\Cx (M)$.
This can be done in two ways.  One possibility is to apply the theorem of Bahturin et al. 
(\cite{BSZ}, Thm. 6) 
on gradings of a matrix algebra. A second, equivalent approach, is to use the classification of cyclotomic 
quantum tori \cite{N, ABFP}. 
We will follow the second approach here, since it can also be adapted for the infinite-dimensional 
(quasifinite) set-up.

 Let us recall the definition of a quantum torus. A {\em quantum torus} is the unital associative
algebra $\R_q = \R_q \left< t_1^\pm, \ldots, t_r^\pm \right>$, whose generators $t_i$ are subject to the 
defining relations
$$ t_i t_j = q_{ij} t_j t_i, \ \ \ i,j = 1, \ldots, r,$$
where $q = (q_{ij})$ is a complex $r \times r$ matrix satisfying $q_{ii} = 1$, $q_{ij} = q_{ji}^{-1}$.   

 We say that the quantum torus $\R_q$ is {\em cyclotomic} if all $q_{ij}$ are complex roots of 1.
 
\medskip
 
\begin{lemma}\label{qtor} 
 Let $\eta_1, \ldots, \eta_r$ be a set of generators of the group $\wH$ of orders
$s_1, \ldots, s_r$, respectively, and let
$T_1, \ldots, T_r$ be invertible operators on $M$ satisfying 
\begin{equation}\label{Ts}
 \rho(\eta_j (a)) = T_j \rho(a) T_j^{-1}, \text{\rm \ \ for all \ } a\in\A,\ j=1,\ldots,r.
\end{equation}
Then the map $t_i\mapsto T_i$ defines a representation of a cyclotomic quantum torus 
$\R_q \left< t_1^\pm, \ldots, t_r^\pm \right>$, with 
$$q_{ij}^{\gcd (s_i, s_j)} = 1.$$
\end{lemma}
\proof
The automorphisms $\eta_i$ and $\eta_j$ commute. Thus for all $a\in \A$,
$$ T_i T_j \rho(a) T_j^{-1} T_i^{-1} = T_j T_i \rho(a) T_i^{-1} T_j^{-1}, $$
or equivalently,
$$ T_i^{-1} T_j^{-1} T_i T_j \rho(a) = \rho(a) T_i^{-1} T_j^{-1} T_i T_j .$$
However, by Schur's Lemma, the centralizer of a finite-dimensional simple module is a multiple of the identity 
map, so
$$ T_i T_j = q_{ij} T_j T_i $$
for some nonzero $q_{ij} \in\C$.

Moreover, the relation $\eta_i^{s_i} = 1$ implies that $T_i^{s_i}$ is a multiple of the identity.
Since $T_i^{s_i} T_j = q_{ij}^{s_i} T_j T_i^{s_i}$, we conclude that
$$q_{ij}^{s_i} = 1 .$$
Likewise, $q_{ij}^{s_j} = 1$, so $q_{ij}^{\gcd (s_i, s_j)} = 1.$
\qed

\bigskip

 We will need elements of the representation theory of cyclotomic quantum tori.
We call an $\R_q$-module $U$ {\it diagonalizable} if every generator $t_i$ is diagonalizable on $U$.
Whenever $\R_q$ is noncommutative, these operators clearly cannot be diagonal with respect to the same basis. 
Let $\CC$ be the category of finite-dimensional diagonalizable $\R_q$-modules.
Since $T_i^{s_i}$ is a multiple of the identity on $M$, we see that $M$ belongs to Category $\CC$.

\begin{lemma}\label{diag}
Let $\R_q$ be a cyclotomic quantum torus, and let $U$ be a module in Category $\CC$.
Then for any $\bm=(m_1,\ldots,m_r) \in \Z^r$, the monomial $\bt^\bm = t_1^{m_1} \ldots t_r^{m_r}$ is 
diagonalizable on $U$.
\end{lemma}
\proof
Choose a positive integer $b$ such that $q_{ij}^b = 1$ for all $i,j = 1, \ldots, r$. Then
the elements $t_i^{\pm b}$ belong to the centre of $\R_q$. Since each of them is diagonalizable on $U$,
they are simultaneously diagonalizable. Moreover, $(\bt^\bm)^b$ belongs to the subalgebra generated
by $t_1^{\pm b}, \ldots, t_r^{\pm b}$, and hence diagonalizable. This implies that $\bt^\bm$ is
diagonalizable as well.
\qed

\medskip 

\begin{proposition}\label{qreps}
 Every module in Category $\CC$ for the cyclotomic quantum torus 
 $\R_q = \R_q \left< t_1^\pm, \ldots, t_r^\pm \right>$ is completely reducible, and all simple $\R_q$-modules 
 in Category $\CC$ have 
the same dimension.  More explicitly, fix an arbitrary simple $\R_q$-module $N$ in $\CC$.  Then all other
simple $\R_q$-modules in $\CC$ can be obtained from $N$ as twists by the rescaling automorphisms of
$\R_q$:
$$ t_j \mapsto \alpha_j t_j, \ \ \ j=1,\ldots,r,$$
where $\alpha_1,\ldots,\alpha_r \in \C^*$.
\end{proposition}
\proof
To prove this proposition, we will put $\R_q$ in a normal form, using the classification of cyclotomic quantum 
tori \cite{N, ABFP}. If we make a change of variables in a quantum torus
using a matrix $\left( a_{ij} \right) \in \GL_r(\Z)$:
$$ \bar t_i = t_1^{a_{i1}} t_2^{a_{i2}} \ldots t_r^{a_{ir}} ,$$
we will get another quantum torus with a new matrix $\bar q_{ij}$. The classification theorem states that 
by means of such a change of variables, a cyclotomic quantum torus can be brought to the normal
form
\begin{equation}\label{normalform}
 \R_q \left<t_1^\pm, \ldots, t_r^\pm \right> \cong \R_{\zeta_1} \left< x_1^\pm, y_1^\pm \right> \otimes \ldots
\otimes \R_{\zeta_\ell} \left< x_\ell^\pm, y_\ell^\pm \right> \otimes \C[z_{2\ell+1}^\pm,\ldots, z_r^\pm] ,
\end{equation}
where the first $\ell$ tensor factors are rank two cyclotomic quantum tori, and the last tensor
factor is a commutative algebra of Laurent polynomials.
By Lemma \ref{diag}, an $\R_q$-module's membership in Category $\CC$ is independent of
the choice of the generators $t_i^\pm$ of $\R_q$.

 The rank two quantum tori $\R_\zeta = \R_\zeta \left< x^\pm, y^\pm \right>$ that appear in (\ref{normalform})
have the following defining relation
\begin{equation*}%\label{xy}
x y = \zeta y x,
\end{equation*}
where $\zeta$ is a primitive complex $d$th root of unity. The centre of $\R_\zeta$ is generated by 
$x^{\pm d}$ and $y^{\pm d}$. Hence the centre of $\R_q$ is the Laurent polynomial algebra:
$$\ZZ (\R_q) = \C [x_1^{\pm d_1}, y_1^{\pm d_1}, \ldots, x_\ell^{\pm d_\ell}, y_\ell^{\pm d_\ell},
z_{2\ell+1}^\pm,\ldots, z_r^\pm ] .$$  

 Let $U$ be a finite-dimensional diagonalizable $\R_q$-module. Since the action of $\ZZ(\R_q)$ on $U$ is
diagonalizable, we may decompose $U$ into the direct sum of $\R_q$-submodules corresponding to various
central characters $\chi: \ZZ(\R_q) \rightarrow \C$:
$$ U = \mathop\bigoplus\limits_{\chi} U^\chi .$$
Note that a central character $\chi$ is determined by its (nonzero) values $\chi(x_i^{d_i})$, $\chi(y_i^{d_i})$,
and $\chi(z_i)$. Hence, by means of the rescaling automorphisms of $\R_q$, we may twist a module with
central character $\chi$ into a module with any other central character. 

We will show that, up to
isomorphism, there exists a unique finite-dimensional simple $\R_q$-module for each central character. 

In order to prove that the module $U$ is is completely reducible, it is enough to consider its single
component $U^\chi$. This component is a module for the quotient $\R_q^\chi$ of $\R_q$ by the central ideal: 
$$\R_q^\chi = \R_q / \big< z - \chi(z)1\ |\ z \in \ZZ(\R_q) \big> .$$
The algebra $\R_q^\chi$ decomposes into the tensor product:
$$\R_q^\chi \cong \R_{\zeta_1}^\chi \left< x_1^\pm, y_1^\pm \right> \otimes \ldots
\otimes \R_{\zeta_\ell}^\chi \left< x_\ell^\pm, y_\ell^\pm \right> .$$
The rank two factors here have the following structure:
$$\R_\zeta^\chi = \R_\zeta \left< x^\pm, y^\pm \right> 
\left/ \left< x^d - \chi(x^d)1,\, y^d - \chi(y^d)1 \right> \right. .$$
In fact, $\R_\zeta^\chi$ is isomorphic to the matrix algebra $M_d(\C)$ under the following isomorphism:

\begin{equation*}
x \longmapsto \big(\chi(x^d)\big)^{1/d} \left(\begin{array}{llllll}
1&0&\cdot&\cdot&0\\
0&\zeta&\cdot&\cdot&0\\
\vdots&\vdots&&\ddots&\vdots\\
0&0&\cdot&\cdot&\zeta^{d-1}
\end{array}\right)
\end{equation*}

\begin{equation}\label{XYmat}
\hbox{\ and\ \ \  \ }
y \longmapsto \left(\begin{array}{llllll}
0&0&\cdot&0&\chi(y^d)\\
1&0&\cdot&0&0\\
0&1&\cdot&0&0\\
\vdots&\vdots&\ddots&\vdots&\vdots\\
0&0&\cdot&1&0
\end{array}
\right).
\end{equation}

This implies that the algebra $\R_q^\chi$ is isomorphic to the matrix algebra of rank $d_1 \cdots
d_\ell$. It is well-known that the matrix algebra is semisimple and has a unique, up to isomorphism,
finite-dimensional simple module of dimension equal to the rank of the matrix algebra.\qed

\bigskip

The following corollary is an immediate consequence of Proposition \ref{qreps}.

\begin{corollary}\label{isotyp} 
Let $N$ be a simple module in Category $\CC$ for the cyclotomic quantum torus 
$\R_q = \R_q \left< t_1^\pm, \ldots, t_r^\pm
\right>$. Then every module $M$ in Category $\CC$ may be written as a finite sum
\begin{equation}\label{NV}
 M \cong N \otimes \left( \mathop\bigoplus\limits_{\alpha\in(\C^*)^r} V^\alpha \right),
\end{equation} 
where the action of $\R_q$ on $N\otimes V^\alpha$ is given by the formula
\begin{equation}\label{NVact}
 t_i (u \otimes v^\alpha) = t_i (u) \otimes \alpha_i v^\alpha, \text{ \ for \ } u\in N, v^\alpha \in V^\alpha.
\end{equation} 
\end{corollary}
\proof
If we take a one-dimensional space $V^\alpha$, then $N \otimes V^\alpha$ is a simple $\R_q$-module,
which is a twist of $N$ by the rescaling automorphism of $\R_q$ given by $\alpha\in(\C^*)^r$.
Thus the decomposition (\ref{NV}) is just an isotypic decomposition of the module $M$.
\qed

\bigskip

We now look at the representation of $\R_q$ described in Lemma \ref{qtor}.
We show that in this module, the elements $\alpha$ appearing in the isotypic decomposition 
of $M$ are not arbitrary, but belong to the group $H$.

\begin{lemma}\label{HNV}
Let $M$ be a finite-dimensional simple module for an $H$-graded algebra $\A$. Suppose that $T_1,\ldots,T_r$ are 
operators on $M$ that satisfy (\ref{Ts}) and define the action of a cyclotomic quantum torus $\R_q$ on $M$. 
Let $N$ be a simple $\R_q$-submodule of $M$. Then the isotypic
decomposition (\ref{NV}) for $M$ can be written as
\begin{equation}\label{NVh}
 M \cong N \otimes \left( \mathop\bigoplus\limits_{h\in H} V^h \right),
\end{equation} 
where the action of $\R_q$ on $N\otimes V^h$ is given by the formula
\begin{equation}\label{NVhact}
 t_i (u \otimes v^h) = t_i (u) \otimes h(\eta_i) v^h, \text{ \ for \ } u\in N, v^h \in V^h.
\end{equation} 
\end{lemma}
\proof
The $H$-grading of the algebra $\A$ yields an $H$-grading of the matrix algebra
\begin{equation}\label{Heigval}
\End_\Cx (M) = \rho(\A) = \mathop\bigoplus\limits_{h\in H} \rho(\A_h),
\end{equation}
where the endomorphisms in $\rho(\A_h)$ satisfy 
\begin{equation}\label{Heigval2}
T_j \rho(a) T_j^{-1} =  \rho(\eta_j (a)) = h(\eta_j) \rho(a) .
\end{equation}
This can be viewed as an eigenspace decomposition of $\End_\Cx (M)$ with respect to conjugations
by  the $T_j$'s. 

 On the other hand, we can construct such an eigenspace decomposition using the isotypic decomposition
(\ref{NV}). In particular, let $N \otimes V^0$ be the isotypic component of $N$ in $M$, and let
$N \otimes V^\alpha$ be another nontrivial component in (\ref{NV}). We claim that
$\Id_N \otimes \Hom_\Cx (V^0, V^\alpha)$ can be viewed as the eigenspace (with eigenvalue $\alpha_j$) for the 
conjugation action of $T_j$ on $\End_\Cx(M)$.  Indeed, extend $S \in \Hom_\Cx (V^0, V^\alpha)$ to
$\End_\Cx (\bigoplus_{\beta} V^\beta)$ by setting $S \left( V^\beta \right) = 0$ when $\beta \neq 0$. 
Then for $n\in N$ and $v\in V^0$, we have 
\begin{eqnarray*}  
T_j (\Id_N \otimes S) T_j^{-1} (n \otimes v) &=&T_j (\Id_N \otimes S) (T_j^{-1} (n) \otimes v)\\
&=&T_j (T_j^{-1} (n) \otimes Sv)\\
&=&\alpha_j n \otimes Sv .
\end{eqnarray*}
Thus $T_j (\Id_N \otimes S) T_j^{-1} = \alpha_j  (\Id_N \otimes S)$.  By (\ref{Heigval}) and (\ref{Heigval2}), 
all eigenvalues of the conjugation action of $T_j$ on $\End_\Cx (M)$ correspond to the elements of $H$, so 
there exists $h\in H$ such that
$\alpha_j = h(\eta_j)$ for all $j = 1, \ldots, r$.
\qed

\bigskip 

Next we explicitly describe the eigenspace decomposition of $\End_\Cx (M)$.
In order to do this, we introduce the group homomorphism
$$ \gamma : \Z^r \rightarrow H,$$
defined by the formula
$$ \gamma({\bf a}) ({\bf \eta}^{\bf b}) = \prod_{i,j = 1}^r q_{ij}^{a_i b_j} , \text{ \ for \ } 
{\bf a},{\bf b} \in \Z^r.$$ 
Here we use the multi-index notation ${\bf \eta^\bb} = \eta_1^{b_1} \ldots \eta_r^{b_r}$ and 
${\bf T^b} = T_1^{b_1} \ldots T_r^{b_r}$.
The fact that $\gamma({\bf a^\prime} + {\bf a^{\prime\prime}}) = \gamma({\bf a^\prime}) + 
\gamma({\bf a^{\prime\prime}})$
is obvious (see (\ref{addnotation}) for the definition of the additive notation), 
so all we need to verify is that this map is well-defined.  That is, we should check that 
$\gamma({\bf a}) (\eta^{{\bf b^\prime}}) = \gamma({\bf a}) (\eta^{{\bf b^{\prime\prime}}})$ whenever 
$ \eta^{{\bf b^{\prime}}} = \eta^{{\bf b^{\prime\prime}}}$ in $H$. 

Taking ${\bf b} = {\bf b^\prime} - {\bf b^{\prime\prime}}$, it is enough to show that 
$\gamma({\bf a}) (\eta^{\bf b}) = 1$ whenever $\eta^{\bf b} = 1$ in $H$.
However,
$$ {\bf T^a T^b} = \gamma({\bf a}) (\eta^{\bf b}) {\bf T^b T^a} .$$
If $\eta^{\bf b} = 1$, then ${\bf T^b}$ is a multiple of the identity and thus commutes with ${\bf T^a}$. 
This implies 
that $\gamma({\bf a}) (\eta^{\bf b}) = 1$.

Note that the kernel of $\gamma$ corresponds to the centre of the quantum torus $\R_q$:
$$\ZZ (\R_q) = \Span \left\{ {\bf t^a}\ |\ {\bf a} \in \ker\, \gamma \right\} .$$
By Schur's Lemma, the centre $\ZZ (\R_q)$ acts on a finite-dimensional simple
$\R_q$-module $N$ by scalar operators.

\ 

 The following proposition is essentially equivalent to the classification of gradings on
the matrix algebra (\cite{BSZ}, Thm. 6). 

\begin{proposition}\label{Hgrad}
Let $M = N \otimes \left( \mathop\bigoplus\limits_{p\in H} V^p \right)$ be the isotypic decomposition
of $M$ as an $\R_q$-module.
Then the space $\hbox{\em End}_\Cx (M)$ decomposes into a direct sum of eigenspaces with respect to conjugation
by the operators $T_j$, $j = 1, \ldots, r$:
$$ \hbox{\em End}_\Cx (M) = \mathop\bigoplus\limits_{h \in H} \rho(\A_h), $$
where 
$$ \rho(\A_h) = \sum_{{\bf a} \in \Z^r} {\bf T^a} \otimes \left( \mathop\bigoplus\limits_{p \in H} 
\Hom_\Cx \left( V^p, V^{p + h + \gamma({\bf a})}\right) \right),$$
is the eigenspace with eigenvalue $h(\eta_j)$ with respect to conjugation by $T_j$.
\end{proposition}
\proof
It is easy to see that the subspaces 
$${\bf T^a} \otimes \left( \mathop\bigoplus\limits_{p \in H} 
\Hom_\Cx \left( V^p, V^{p + h + \gamma({\bf a})}\right)
\right)$$
span $\End_\Cx (M) = \End_\Cx(N) \otimes \End_\Cx \left(\bigoplus_{h\in H} V^h \right)$. It remains only to 
verify that these subspaces are the eigenspaces for the conjugation by $T_j$'s. 

Let $S \in \Hom_\Cx \left( V^p, V^{p + h + \gamma({\bf a})} \right)$, $n\in N$, $v\in V^p$. Then
\begin{eqnarray*}
T_j ({\bf T^a} \otimes S) T_j^{-1} (n \otimes v) &=&
T_j ({\bf T^a} \otimes S) (T_j^{-1} (n) \otimes p^{-1} (\eta_j) v)\\
&=& 
T_j ({\bf T^a} T_j^{-1} (n) \otimes p^{-1} (\eta_j) Sv)\\
 &=& p(\eta_j) h(\eta_j) \gamma({\bf a}) (\eta_j) p^{-1} (\eta_j) T_j {\bf T^a} T_j^{-1} (n) \otimes Sv\\
&=& h(\eta_j) ({\bf T^a} \otimes S) (n \otimes v) .
\end{eqnarray*}
\qed

\bigskip

Using this result, we can now describe the thin $H$-coverings of the module $M$.
 
\begin{theorem}\label{Hcov}
 Let $M$ be a finite-dimensional simple $\A$-module, and let $\wH$ be a finite abelian subgroup 
of $\Aut\,\A$, such that the twisted modules $M^\eta$ are isomorphic to $M$ for all
$\eta\in\wH$. 
Let 
$$ M  = \mathop\bigoplus\limits_{h\in H} \big(N \otimes V^h\big) $$
be the isotypic decomposition of $M$ with respect to the action of the cyclotomic quantum torus $\R_q$ given by 
(\ref{NVhact}).
Then up to equivalence, thin $H$-coverings of $M$ are parametrized by one-dimensional subspaces in $N$. 
Namely, let $n\in N$, $n \neq 0$. Define
\begin{equation}\label{Mh}
M_h = \sum_{{\bf a}\in\Z^r} {\bf T^a} (n) \otimes V^{h + \gamma({\bf a})} .
\end{equation}
Then $\left\{ M_h\ |\ h\in H \right\}$ is a thin $H$-covering of $M$, and every thin $H$-covering 
of $M$ is equivalent to a covering of this type.
\end{theorem}

\begin{remark}{\em 
The sum in (\ref{Mh}) may be replaced with a finite sum taken over ${\bf a}\in \Z^r / \ker \gamma$.}
\end{remark}

\noindent
\proof 
By Proposition \ref{Hgrad}, 
$$ \rho(\A_m) = \sum_{{\bf b}\in\Z^r} {\bf T^b} \otimes \left( \mathop\bigoplus\limits_{p\in H} \Hom_\Cx 
\big(V^p, V^{p+m+\gamma({\bf b})} \big) \right)$$
for any $m\in H$.  Thus 
\begin{eqnarray}\label{dom}
\rho(\A_m) M_h &=&\sum_{{\bf a},{\bf b}\in\Z^r} {\bf T^b T^a} (n) \otimes V^{h+m+\gamma({\bf a})+
\gamma({\bf b})}\\
&=&\sum_{{\bf c}\in\Z^r} {\bf T^c} (n) \otimes V^{h+m+\gamma({\bf c})}\\
&=& M_{h+m}.
\end{eqnarray}

The subspaces $M_h$ clearly span $M$:
\begin{eqnarray*}
\sum_{h\in H} M_h &=& \sum_{h\in H} \sum_{{\bf a}\in\Z^r} {\bf T^a} (n) \otimes V^{h + \gamma({\bf a})}\\
&=& \left( \sum_{{\bf a}\in\Z^r} {\bf T^a}(n) \right) \otimes 
\left( \mathop\bigoplus\limits_{h\in H} V^h \right)\\
&=& N \otimes  \left( \mathop\bigoplus\limits_{h\in H} V^h \right)\\
 &=& M .
\end{eqnarray*}
Thus
$\left\{ M_h\ |\ h\in H \right\}$ is an $H$-covering of $M$.

We will use Lemma \ref{obvious} to show that this covering is thin.  All that remains is to determine the 
simple $\rho(\A_0)$-submodules of $M$.  Let $L \subseteq M$ be a simple $\rho(\A_0)$-submodule, and let $u$ be 
a nonzero vector in $L$.  Fix a basis $\{ v_{hj}\ |\ 1\leq j\leq \dim V^h \}$ for each of the spaces $V^h$.  
Expand $u$ according to this basis:
$$ u = \sum\limits_{h,j} n_{hj} \otimes v_{hj}\hbox{\  with\ } n_{hj} \in N.$$ 
Since the algebra $\rho(\A_0)$ contains a subalgebra
$$\Id_N \otimes \left( \mathop\bigoplus\limits_{h\in H} \End_\Cx \left(V^h \right) \right),$$
we see that $L$ contains $n_{hj} \otimes v_{hj}$ for all $h,j$.
Thus without loss of generality, we may assume that $u = n \otimes v$, where
$n \in N$, $v\in V^h$. Since
$$\rho(\A_0) = \sum_{{\bf a}\in\Z^r} {\bf T^a} \otimes \left( \mathop\bigoplus\limits_{p\in H} 
\Hom_\Cx \left( V^p, V^{p+\gamma({\bf a})} \right) \right),$$
$L$ contains the subspace
\begin{equation}\label{simpleA}
U^h (n) = \sum_{{\bf a}\in\Z^r} {\bf T^a} (n) \otimes V^{h + \gamma({\bf a})} .
\end{equation}
The subspace $U^h (n)$ is $\rho(\A_0)$-invariant, and hence every simple $\rho(\A_0)$-sub-
module in $M$
coincides with $U^h (n)$ for some $n\in N$, $h\in H$, $n\neq 0$, $V^h \neq 0$. 

Let us show that $U^h (n)$ is, in fact, a 
simple $\rho(\A_0)$-module for all $n\in N$, $h\in H$ such that $n \neq 0$ and $V^h \neq 0$.
Fix $\bb\in\Z^r$. The intersection of
$U^h(n)$ with $N \otimes V^{h+\gamma(\bb)}$ is 
$$\sum\limits_{\ba\in\ker\gamma} \T^\bb \T^\ba (n) \otimes V^{h+\gamma(\bb)}.$$ 
The centre $\ZZ (\R_q)$ acts on 
the simple $\R_q$-module $N$ by scalars.  Thus for all $\ba\in\ker\gamma$, $\T^\ba(n)$ is a multiple of $n$, 
so the intersection of $U^h(n)$ with $N \otimes V^{h+\gamma(\bb)}$ is  
$\T^\bb(n) \otimes V^{h+\gamma(\bb)}$. By the argument given above, every nonzero vector in $U^h(n)$
generates $U^h(n)$. Therefore the spaces $U^h(n)$ with $V^h \neq 0$, $n \neq 0$, exhaust all the simple 
$\rho(\A_0)$-submodules in $M$.

 This proves that the components $M_h$ of the $H$-covering (\ref{Mh}) are simple $\rho(\A_0)$-modules 
whenever $M_h \neq 0$.  Thus for any nonzero $m\in M_h$, we have 
%\begin{eqnarray*}
$$\rho(\A_{g-h})m\supseteq \rho(\A_{g-h})\rho(\A_0)m=\rho(\A_{g-h})M_h=M_g,$$
using (\ref{dom}).  Hence the covering $\{M_h\ |\ h\in H\}$ is thin by Lemma \ref{obvious}.

Since the described components of the thin coverings exhaust all the simple $\rho(\A_0)$-submodules
in $M$, we obtain that up to equivalence, these are all the thin $H$-coverings of $M$.
\qed

\bigskip

\begin{corollary}\label{griso}
 Let $M$ be a finite-dimensional simple $\A$-module, and let $\wH$ be a finite abelian subgroup  of $\Aut\,\A$, such that the twisted modules $M^\eta$ are isomorphic to $M$ for all
$\eta\in\wH$. 
Let $\{ M_h\ |\ h\in H\}$ and $\{ M^\prime_h\ |\ h\in H\}$ be two thin coverings of $M$. Then for some
$g \in H$, the graded-simple modules
$$\widetilde{M} = \mathop\bigoplus\limits_{h\in H} M_{h} \text{\rm \ and \ } 
\widetilde{M}^\prime = \mathop\bigoplus\limits_{h\in H} M^\prime_{h+g}$$
are isomorphic as graded $\A$-modules. 
\end{corollary}
\proof
Clearly, we are allowed to replace the thin coverings $\{M_h\}$ and $\{M^\prime_h\}$ with any equivalent thin coverings.  By Theorem \ref{Hcov}, we may assume that for all $h\in H$,
$$\widetilde{M}_{h} = \sum\limits_{\ba\in\Z^r} \T^\ba (n) \otimes V^{h+\gamma(\ba)} \text{\rm \ and \ } 
\widetilde{M}^\prime_{h} = \sum\limits_{\ba\in\Z^r} \T^\ba (n^\prime) \otimes V^{h+\gamma(\ba)}$$
for some nonzero $n, n^\prime \in N$.  %Let $g=p^\prime-p$.

We construct a grading-preserving isomorphism $\psi: \widetilde{M} \rightarrow  \widetilde{M}^\prime$
defined on $\widetilde{M}_{h}$ by
$$\psi \left( \T^\ba (n) \otimes v \right) = \T^\ba (n^\prime) \otimes v $$
for all $h\in H, \ba \in\Z^r, v\in V^{h+\gamma(\ba)}$. Using Proposition \ref{Hgrad}, it is easy to see that
$\psi$ commutes with the action of $\A$. 
\qed

\medskip

We now consider some basic examples of thin $H$-coverings.
 
\begin{example} \label{cycex} {\em Suppose $\A$ is graded by a finite cyclic group 
$\wH=\langle \eta\rangle\subseteq\Aut\A$, and suppose that there is a (normalized) isomorphism 
$T:\ M\rightarrow M^\eta$.  Then, up to equivalence, there is a unique thin $H=\Hom(\wH,\C^*)$-covering 
$\{M_h\}$ of $M$, where $M_h$ is the eigenspace
$$M_h=\{u\in M\ |\ Tu=h(\eta)u\}.$$
In this case, the module $M$ admits an $H$-grading: $M = \mathop\bigoplus\limits_{h\in H} M_h$.}
\end{example}
\proof
Since the sum $M=\sum_{h\in H}M_h$ is direct, the thinness of the covering $\{M_h\}$ in Example \ref{cycex} is 
obvious.  However, even in this simplest of examples, the {\em uniqueness} of the thin covering requires 
Theorem \ref{Hcov}.
 Note that, in this case, the quantum torus $\R_q$ is just the algebra of Laurent polynomials in a single 
variable, so $\dim N=1$, and $M$ may be identified with $\bigoplus_{h\in H}V_h$.\qed

\bigskip

From the statement of Theorem \ref{Hcov}, it is clear that uniqueness of thin coverings is the exception, 
rather than the rule, as is illustrated by the next pair of examples.  %Thin coverings also provide a finer 
%invariant than graded isomorphism.

\bigskip

\begin{example} \label{matrixex} {\em Let $\A=M_n(\K)$ be the associative algebra of $n\times n$ matrices.
If $\A=\bigoplus_{g\in G}\A_g$ is any grading satisfying conditions (i) and (ii) of Theorem 
\ref{unicite}, then by Theorem \ref{completeness}, every graded-simple module for $\A$ comes from a 
thin covering of a simple module for $\A$.  By the Artin-Wedderburn Theorem, there is only one simple 
module for $\A$; and by Theorem \ref{unicite}, there is only one graded-simple module 
(up to graded-isomorphism) that can come from the thin coverings of this module.  Therefore (up to graded 
isomorphism), $\A$ has only one graded-simple module.

In spite of the uniqueness of the graded-simple module for $\A$, the thin coverings for $\A$ are far 
from unique.  For example, $\A$ has a $\Z_n\times\Z_n$-grading with
$$\A_{(\overline{a},\overline{b})}:=\hbox{Span}\{E^aF^b\}$$
where 
$$E:= \left(\begin{array}{llllll}
1&0&\cdot&\cdot&0\\
0&\zeta&\cdot&\cdot&0\\
\vdots&\vdots&&\ddots&\vdots\\
0&0&\cdot&\cdot&\zeta^{d-1}
\end{array}\right),
%\hbox{\ and\ }
\hbox{\quad\ }
F:= \left(\begin{array}{llllll}
0&0&\cdot&0&1\\
1&0&\cdot&0&0\\
0&1&\cdot&0&0\\
\vdots&\vdots&\ddots&\vdots&\vdots\\
0&0&\cdot&1&0
\end{array}
\right),
$$
$\zeta$ is a primitive $n$th root of unity, and $\ol{a}$ (resp., $\ol{b}$) is the image of $a\in\Z$ (resp., 
$b\in\Z$) under the canonical homomorphism $\Z\rightarrow\Z_n=\Z/n\Z$.

Then the graded ring $\A$ satisfies conditions (i) and (ii) of Theorem \ref{unicite}, so any thin covering 
of the natural module $M=\K^n$ induces the same (unique) graded-simple module of $\A$ (relative to the 
grading group $G=\Z_n\times\Z_n$).  Let $v$ be a nonzero vector in $M$.  Then since $M$ is irreducible, 
$$M_g:=\A_g v$$
defines a covering of $M$.  Let $\{N_g\ |\ g\in G\}$ be another covering of $M$ with $N_g\subseteq M_g$ for 
all $g\in G$.  Then $N_h\neq 0$ for some $h\in G$, so $N_g=\A_{g-h}.N_h\neq 0$ for all $g\in G$, since the 
one-dimensional space $\A_{g-h}$ is spanned by the invertible matrix $E^aF^b$ (where $g-h=(\ol{a},\ol{b})$).  
Moreover, $N_g$ is a nonzero subspace of the one-dimensional space $M_g$, so $N_g=M_g$ for all $g\in G$.  
Thus the covering $\{M_g\}$ is thin.  

Any choice of $v$ results in a thin covering consisting of $n^2$ one-dimensional 
subspaces $M_g$ of $M$.  The space $M_0$ is the span of the vector $v$.  Since there are infinitely many distinct
one-dimensional subspaces we could choose for $M_0$, there are infinitely many nonequivalent thin 
coverings of $M$.}\end{example}

More generally, using Theorem \ref{Hcov}, we have the following example.

\begin{example}{\em 
Suppose $\wH$ is generated by two automorphisms $\eta_1, \eta_2$ of order $d$.
Let $T_1, T_2$ be isomorphisms from $M$ to $M^{\eta_1}$, $M^{\eta_2}$, respectively.
Assume that these operators satisfy the relation
$$T_1 T_2 = \zeta T_2 T_1,$$
where $\zeta$ is a primitive root of unity of order $d$.  Consider the isotypic decomposition
$$ M = N \otimes \left( \mathop\bigoplus\limits_{h \in H} V^h \right) $$
of $M$ as an $\R_\zeta$-module. We obtain
from Theorem \ref{Hcov} that the thin $H$-coverings are parametrized by one-dimensional subspaces in 
a $d$-dimensional space $N$. For any nonzero vector $n\in N$ we get a thin $H$-covering with
$$M_h = \sum_{{\bf a}\in\Z^2} {\bf T^a} (n) \otimes V^{h + \gamma({\bf a})} .$$
Up to equivalence, these are all thin $H$-coverings of $M$.}
\end{example}

\medskip

We are finally ready to describe the thin $G$-coverings of $M$ in the general case--that is, when the twists of 
$M$ by
automorphisms in $\wG$ are not necessarily isomorphic to $M$. The following theorem shows that $G$-coverings
of $M$ can be described via $H$-coverings, where $\wH$ is the subgroup of $\wG$ defined in Lemma 
\ref{Hsubgroup}. In order to distinguish between $G$-coverings and $H$-coverings (resp. gradings), 
we will add an appropriate superscript to our notation.

 The inclusion $\wH \subseteq \wG$ yields a natural epimorphism
$$ \psi: \ \Hom(\wG,\C^*) \rightarrow \Hom(\wH,\C^*) $$
by restriction of the maps to $\wH$. 

\medskip

\begin{theorem}\label{Gcov}
Let $M$ be a finite-dimensional simple $\A$-module, and let $\wG$ be a finite abelian
subgroup of $\Aut\,\A$. Let $\wH$ be its subgroup defined in Lemma \ref{Hsubgroup}. Then

(i)\quad\ $\rho(\A^G_g) = \rho\big(\A^H_{\psi (g)}\big)$ for all $g\in G$.

(ii)\quad\!\!\! There is a bijective correspondence between thin $G$-coverings of $M$ and thin $H$-coverings
of $M$. Given a thin $H$-covering $\{ M^H_h \}$ we get a thin $G$-covering
$\{ M^G_g \}$ with $M^G_g = M^H_{\psi(g)}$. Every thin $G$-covering is of this form.

(iii)\quad\!\!\!\!\! All graded-simple $\A$-modules associated with various thin coverings of $M$ (as in 
Theorem \ref{grmodlm})
are isomorphic (up to shifts in the gradings on the modules).
\end{theorem}
\proof Let $\tau_1,\ldots,\tau_\ell$ be representatives of the cosets of $\wH$ in $\wG$, where 
$\ell = |\wG|/|\wH|$.
By Lemma \ref{Hsubgroup}, the modules $M^{\tau_1}, \ldots, M^{\tau_\ell}$ are pairwise
non-isomorphic.
Consider the representation $\hro$ of $\A$ on 
$$ \hM =  M^{\tau_1} \oplus  \cdots \oplus M^{\tau_\ell } .$$
An element $a\in \A^G_g$ acts on $\hM$ by
$$\hro (a) = \left( g(\tau_1) \rho(a),  \ldots, g(\tau_\ell )\rho(a) \right) .$$
Thus $\dim \hro(\A^G_g) = \dim \rho(\A^G_g)$.  Moreover, by the definition of the map $\psi$, 
$\A^G_g \subseteq \A^H_{\psi(g)}$, so $\rho(\A^G_g) \subseteq \rho\big(\A^H_{\psi(g)}\big)$.  Hence
\begin{eqnarray}\label{dims}
\dim \hro(\A) &\leq& \sum_{g\in G} \dim \hro(\A^G_g)\\
&=& \sum_{g\in G} \dim \rho(A^G_g)\\ 
&\leq& \sum_{g\in G} \dim \rho(\A^H_{\psi(g)})\\
&=& \big(|G|/|H|\big) \sum_{h\in H} \dim \rho(\A^H_h)\\
&\leq& \ell  \dim \End_\Cx (M) .\label{dimf}
\end{eqnarray}
However, by the density theorem,
$$\hro (\A) = \End_\Cx(M^{\tau_1}) \oplus \cdots \oplus \End_\Cx (M^{\tau_\ell }) .$$
This implies that all inequalities in (\ref{dims})--(\ref{dimf}) are in fact equalities, and in particular,
$$ \dim \rho(\A^G_g) = \dim \rho (\A^H_{\psi(g)}) .$$  
This completes the proof of Part (i) of the theorem.

 To prove the second part, we note that $\rho(\A^G_0) = \rho(\A^H_0)$. This means that the simple 
$\rho(\A^G_0)$-submodules in $M$ are the simple $\rho(\A^H_0)$-submodules.
Let $\{ M^G_g \}$ be a thin $G$-covering of $M$. Without loss of generality, we may
assume that $M^G_0 \neq 0$. Then by Lemma \ref{obvious}, $M^G_0$ is a simple $\rho(\A^H_0)$-submodule.
We have seen in the proof of Theorem \ref{Hcov} that every simple $\rho(\A^H_0)$-submodule in $M$
occurs as a component of a thin $H$-covering.  Since we can always replace such a covering by an equivalent 
thin covering, we see that every simple $\rho(\A_0^H)$-submodule in $M$ occurs as the $0$-component of a thin 
$H$-covering.  Thus there exists a thin $H$-covering
$\{M^H_h \}$, such that $M^H_0 = M^G_0$. Moreover, $M^H_h = \rho\left( \A^H_h \right) M^H_0$ for all $h \in H$.
By Lemma \ref{obvious} and Part (i), we get that
$$M^G_g = \rho\left( \A^G_g \right) M^G_0 = \rho\left( \A^H_{\psi(g)} \right) M^H_0 = M^H_{\psi(g)} .$$
Thus every thin $G$-covering comes from a thin $H$-covering. Conversely, 
if we construct a $G$-covering from a thin $H$-covering by taking $M^G_g = M^H_{\psi(g)}$,
then Lemma \ref{obvious} implies that this $G$-covering is thin.

Part (iii) follows from (i), (ii), and Corollary \ref{griso}.
\qed 

\bigskip

When $\wG$ is finite cyclic, Theorem \ref{Gcov} reduces to the following eigenspace decomposition.
\begin{example} {\em Suppose $\wG$ is a finite cyclic group generated by $\eta\in\Aut\A$, $\wH\subseteq\wG$ is 
generated by $\eta^s$, and $T:\ M\rightarrow M^{\eta^s}$ is a (normalized) module isomorphism.  Then 
(up to equivalence) $M$ has a unique thin $G$-covering given by 
$$M_g=\{u\in M\ |\ Tu=g(\eta)^su\}.$$}
\end{example}

\section{Thin coverings of quasifinite modules}

In this section we generalize our results to an infinite-dimensional setting.
We will assume that the algebra $\A$ has an additional grading by an abelian group $Z$,
$$\A = \mathop\bigoplus\limits_{k\in Z} \A_k,$$
and that $M$ is a $Z$-graded $\A$-module,
$$M = \mathop\bigoplus\limits_{n\in Z} M_n,$$
with finite-dimensional graded components $M_n$.  Such modules are called {\em quasifinite}.

Note that $M$ may be infinite-dimensional, since we allow $Z$ to be infinite.
We are primarily interested in the case when $Z = \Z$, but our results hold for arbitrary
abelian grading groups.

We will work in the categories of $Z$-graded algebras and modules and consider only those automorphisms 
$\sigma$ which preserve the $Z$-grading.  That is, 
$$\sigma(\A_k) = \A_k.$$

Let $M$ be a quasifinite graded-simple module (with no proper $Z$-graded submodules).  In this section, we 
show that Theorems \ref{Hcov} and \ref{Gcov} still hold in this more general setting.  We start with some 
obvious adjustments to the definition of a thin $G$-covering.

Let $\wG$ be a finite abelian group of ($Z$-grading-preserving) automorphisms of $\A$. Then
$\A$ is graded by the group $G = \Hom (\wG, \C^*)$, and the two gradings are compatible:
$$ \A = \mathop\bigoplus\limits_{(k,g)\in Z \times G} \A_{k,g} .$$
A {\em covering} of $M$ is a set of $Z$-graded subspaces $\{ M_g\ |\ g\in G \}$ of $M$ 
that span $M$ and satisfy
\begin{equation}\label{quasiAM}
\rho(\A_g) M_h \subseteq M_{g+h} , \ \ \text{\rm for all\ } g,h \in G.
\end{equation}
As before, the set of coverings is partially ordered, and the minimal elements of the poset of coverings
are called {\em thin coverings}.

\medskip

For all $n,k\in Z$, let $\rho_{nk}$ be the restriction of the action of $\A_{k-n}$ to $M_n$:
$$ \rho_{nk} : \A_{k-n} \rightarrow \Hom_\Cx(M_n, M_k) .$$
It is easy to see that $G$-coverings can be defined in terms of $\rho_{nk}$. Indeed, condition
(\ref{quasiAM}) can be replaced by
$$\rho_{nk} (\A_{k-n,h}) M_{n,g} \subseteq M_{k,g+h},$$
and ``local''
information about $\rho_{nk} (\A_{k-n}) \subseteq \Hom_\Cx(M_n, M_k)$ can substitute for ``global'' 
information about $\rho(\A)\subseteq \End_\Cx(M).$  The proofs of the previous section are essentially based
on Schur's Lemma and the Jacobson Density Theorem.  We will use the quasifinite version of the density 
theorem presented in the appendix.

\begin{remark}\label{qmax}{\em 
The assumption on the existence of a maximal submodule made in Theorem \ref{converse}
holds for the quasifinite modules described above.}
\end{remark}
\proof
 Let $M$ be a $Z$-graded quasifinite module, which also has
a compatible $G$-grading and is graded-simple (has no proper $Z \times G$-homogeneous submodules).
Fix $n\in Z$ such that $M_n \neq 0$.
Note that by Zorn's Lemma, $M$ has a maximal $Z$-graded submodule.
Indeed, if $U$ is a proper $Z$-graded submodule in $M$ then $U_n \neq M_n$.
Otherwise, if $U_n = M_n$ then the space $U_n$ is homogeneous both with respect to $Z$ and $G$ and will
generate the whole module $M$ since it is graded-simple.

Now if we consider an increasing chain of $Z$-graded submodules in $M$,
$$ \cdots \subseteq U^{(i)} \subseteq U^{(i+1)} \subseteq \cdots ,$$ 
the corresponding chain of components 
$$ \cdots \subseteq U^{(i)}_n \subseteq U^{(i+1)}_n \subseteq \cdots  \subseteq M_n$$ 
will stabilize, since $\dim M_n < \infty$.
Thus the union $\mathop\bigcup\limits_{i} U^{(i)}$ is a $Z$-graded submodule in $M$, with $n$th component
properly contained in $M_n$. Applying Zorn's Lemma, we conclude that the partially ordered set of $Z$-graded 
proper submodules of $M$ contains a maximal element.
\qed

\bigskip

Our main results about quasifinite modules, Theorems \ref{quasiHcov} and \ref{quasiGcov}, can be proven 
using the same arguments as for the finite case.  Only the proof of Theorem \ref{quasiGcov}(i) requires 
nontrivial modification, since the corresponding finite result, Theorem \ref{Gcov}(i), used a dimension argument. 
We write down the details of the modified proof below, leaving the other proofs as straightforward 
exercises to the reader.
     
\begin{theorem}\label{quasiHcov}
 Let $M$ be a $Z$-graded-simple quasifinite module for a $Z$-graded algebra $\A$. 
Let $\wH$ be a finite abelian group 
of automorphisms of $\A$ preserving the $Z$-grading, such that 
there exist $Z$-grading-preserving isomorphisms $T_\eta$ between $M$ and $M^\eta$ for all
$\eta\in\wH$. 
Let 
$$M  = \mathop\bigoplus\limits_{h\in H} N \otimes V^h$$
be the isotypic decomposition of $M$ with respect to the action of the cyclotomic quantum torus 
$\R_q$ given by the operators $T_\eta$. 
In this decomposition, each subspace $V^h$ is $Z$-graded.
Then up to equivalence, thin $H$-coverings of $M$ are parametrized by one-dimensional
subspaces in $N$. 
Namely, for $n\in N$, $n \neq 0$, define 
$$M_h = \sum_{a\in\Z^r} {\bf T^a} (n) \otimes V^{h + \gamma(a)} .$$
Then $\left\{ M_h\ |\ h\in H \right\}$ is a thin $H$-covering of $M$, and every thin $H$-covering 
of $M$ is equivalent to a covering of this type.\qed
\end{theorem}

\smallskip

\begin{theorem}\label{quasiGcov}
Let $M$ be a $Z$-graded-simple quasifinite module for a $Z$-graded algebra $\A$.
Let $\wG$ be a finite abelian
group of automorphisms of $\A$ preserving the $Z$-grading. 
Let $\wH$ be its subgroup defined in Lemma \ref{Hsubgroup}. Then

(i)\quad\ $\rho_{nk}( \A^G_{k-n,g}) = \rho_{nk}(\A^H_{k-n,\psi (g)})$ for all $g\in G$.

(ii)\quad\!\!\! There is a bijective correspondence between thin $G$-coverings of $M$ and thin $H$-coverings
of $M$. Given a thin $H$-covering $\{M^H_h \}$ we get a thin $G$-covering
$\{ M^G_g \}$ with $M^G_g := M^H_{\psi(g)}$. Every thin $G$-covering of $M$ is of this form.

(iii)\quad\!\!\!\!\! All graded-simple $\A$-modules associated with various thin coverings of $M$ (as in 
Theorem \ref{grmodlm})
are isomorphic (up to shifts in the $G$-gradings on the modules).
\end{theorem}
\proof
We will give a proof for Part (i) of the theorem.  Parts (ii) and (iii) may be proven analogously to 
Theorem \ref{Gcov}. 
% It follows from the definition of the map $\psi$ that $\A_g \subseteq \A^H_{\psi(g)}$,
% and so $\rho(\A_g) \subseteq \rho(\A^H_{\psi(g)})$.

Let $\tau_1,\ldots,\tau_\ell $ be representatives of the cosets of $\wH$ in $\wG$, where $\ell  = |\wG|/|\wH|$.
The modules $M^{\tau_1}, \ldots, M^{\tau_\ell }$ do not admit ($Z$-grading-preserving) isomorphisms
between any distinct pair of them.

Consider the representation $\hro$ of $\A$ on 
$$ \hM =  M^{\tau_1} \oplus  \cdots \oplus M^{\tau_\ell } .$$
For $n,k\in Z$, let $\hro_{nk}$ be the restriction of $\hro(\A_{k-n})$ to 
$\hM_n=M_n^{\tau_1} \oplus  \cdots \oplus M_n^{\tau_\ell }$:
$$\hro_{nk}: \A_{k-n} \rightarrow \Hom_\Cx( \hM_n , \hM_k ) .$$
An element $a\in \A^G_{k-n,g}$ acts on $\hM_n$ by
$$\hro_{nk} (a) = \big( g(\tau_1) \rho_{nk}(a),  \ldots, g(\tau_k)\rho_{nk}(a) \big) .$$
Thus $\dim \hro_{nk}(\A^G_{k-n,g}) = \dim \rho_{nk}(\A^G_{k-n,g})$. Hence
$$ \dim \hro_{nk}(\A_{k-n}) \leq \sum_{g\in G} \dim \hro_{nk}(\A^G_{k-n,g}) = \sum_{g\in G} \dim
\rho_{nk}(\A^G_{k-n,g})$$
\begin{equation}\label{quasidims}
 \leq \sum_{g\in G} \dim \rho_{nk}\big(\A^H_{k-n,\psi(g)}\big)=\big(|G|/|H|\big)
 \sum_{h\in H}\dim\rho_{nk}\big(\A_{k-n,h}^H\big),
\end{equation}
using the fact that $\A_{k-n,g}^G\subseteq\A_{k-n,\psi(g)}^H$ by the definition of $\psi$.  The space 
$\rho_{nk}(\A_{k-n})$ decomposes into a direct sum
$$\rho_{nk}(\A_{k-n}) = \mathop\bigoplus\limits_{h\in H} \rho_{nk}(\A^H_{k-n,h})$$
(cf.~Proposition \ref{Hgrad}), so 
$$\big(|G|/|H|\big)\sum_{h\in H}\dim\rho_{nk}\big(\A_{n-k,h}^H\big)=\ell\dim\rho_{nk}\big(\A_{k-n}\big).$$
But by the Quasifinite Density Theorem (\ref{qJDT}), 
$$ \rho_{nk} (\A_{k-n}) = \Hom_\Cx(M_n, M_k) $$
and
$$ \hro_{nk} (\A_{k-n}) = \mathop\bigoplus\limits_{j=1}^\ell \Hom_\Cx(M^{\tau_j}_n, M^{\tau_j}_k),$$
so $\dim\hro_{nk}(\A_{n-k})=\ell\dim\rho_{nk}(\A_{n-k}).$

Therefore all the inequalities in (\ref{quasidims}) are in fact equalities, and in particular,
$$ \dim \rho_{nk}(\A^G_{k-n,g}) = \dim \rho_{nk} (\A^H_{k-n,\psi(g)}) .$$  
This completes the proof of Part (i).
\qed

\section{Application to multiloop Lie algebras}

Let $\wG$ be a finite abelian group of grading-preserving automorphisms of 
a $Z$-graded complex Lie algebra  $\g$. Consider the dual group $G = \Hom (\wG, \C^*)$.
The Lie algebra $\g$ is then graded by the group $Z \times G$.
By choosing a set of $N$ generators of $G$, we get an epimorphism $\Z^N \rightarrow G$.
We will denote by $\ol{\ba}$ the image of an element $\ba \in \Z^N$ under this map.

Consider the multiloop Lie 
algebra $L(\g,\wG)$: 
$$L(\g,\wG)=\bigoplus_{\ba\in\Z^{N}}\g_{\ol{\ba}}\ot {\bf t}^\ba,$$
which is a Lie subalgebra of $\g \ot \C [t_1^\pm, \ldots, t_N^\pm]$.

Let $M$ be a quasifinite $Z$-graded-simple module for $\g$. Using Theorem \ref{quasiGcov},
we can construct a $G$-covering $\{ M_g \}$ of $M$. Then the space
$$\bigoplus_{\ba\in\Z^{N}} M_{\ol{\ba}}\ot {\bf t}^\ba$$  
becomes a $(Z \times \Z^N)$-graded-simple module for the multiloop algebra $L(\g,\wG)$.

\appendix
\section{Appendix}

In this appendix, we prove the quasifinite density theorem used in the proof of Theorem
\ref{quasiGcov}.

 Let $M$ be a $Z$-graded module over a $Z$-graded associative ring $R$. 
 We call such a module {\em graded-semisimple} if it is a direct sum of a finite
number of graded-simple modules. Consider the set of endomorphisms of $M$ that shift the grading by an 
element $p\in Z$:
$$ \End (M)_p = \left\{ X \in \End(M)\ |\ X M_n \subseteq M_{n+p} \text{\ for all \ } n \in Z \right\} .$$
In this notation, $\End (M)_0$ is the ring of grading-preserving endomorphisms of $M$.

 The following is a graded version of the usual density theorem of Jacobson and 
Chevalley (cf.~\cite{lam}, Thm. 4.11.16, for instance):

\begin{theorem}[Graded Density Theorem]\label{JDT}Let $M$ be a  
$Z$-graded-semi-
simple module over a $Z$-graded associative ring $R$. 
Let $K = \End_R (M)_0$ and let $E= \mathop\bigoplus\limits_{p\in Z} \End_K(M)_p$.  
Suppose that $\{x_1,\ldots,x_m\} \subseteq M_s$ for some $s\in Z$.
Then for every $f \in E$ there exists an element $r\in R$ 
such that $rx_j=f(x_j)\hbox{\ for\ }j=1,2,\ldots,m.$
\end{theorem}

\noindent
\proof The argument given in \cite{lam} is valid in the graded set-up as well.
\qed

\begin{theorem}[Quasifinite Density Theorem]\label{qJDT}
Let $\A$ be an associative algebra graded by an abelian group $Z$.
Let $M^{(1)},...,M^{(n)}$ be $Z$-graded-simple quasifinite $\A$-modules, such that
there is no grading-preserving isomorphism between any distinct pair of them. Fix $s, p \in Z$, and for each 
$i=1,\ldots,n$, let $\{x^i_j \}_{j=1,\ldots,m_i}$ be a $\Cx$-linearly independent subset of $M^{(i)}_s$, 
and let $\{y^i_j \}_{j=1,\ldots,m_i}$ be a set of vectors in $M^{(i)}_{s+p}$.
Then there exists an element $a \in \A_p$, such that $ax^i_j = y^i_j$ for all $i,j$.
\end{theorem}
\proof
By Schur's Lemma, every nonzero grading-preserving homomorphism 
of graded-simple modules is an isomorphism.  Thus
$$\Hom_\A \left( M^{(i)}, M^{(\ell)} \right)_0 = 0\hbox{\  for\ } i \neq \ell.$$
Applying Schur's Lemma again,
we see that every grading-preserving endomorphism of a quasifinite graded-simple module is a multiple
of the identity. Take $M= M^{(1)} \oplus \cdots \oplus M^{(n)}$. Then
$$ K = \End_\A (M)_0 = \C \pi_1 \oplus \cdots \oplus \C \pi_n ,$$
where $\pi_j$ is the projection $\pi_j : M \rightarrow M^{(j)}$. It follows that 
$$ \End_K (M)_p = \mathop\bigoplus_{j=1}^n \End_\C \big(M^{(j)}\big)_p ,$$
and now Theorem \ref{qJDT} follows immediately from Theorem \ref{JDT}.
\qed

%\bigskip

%\small
%\noindent
%Carleton University,
%School of Mathematics and Statistics,
%1125 Colonel By Drive,
%Ottawa, Ontario, Canada K1S 5B6\\
%Email: billig@math.carleton.ca
%\bigskip

%\noindent
%University of Ottawa,
%Department of Mathematics and Statistics,
%585 King Edward Ave.,
%Ottawa, Ontario, Canada K1N 6N5\\
%Email: mlau@uottawa.ca
\end{document}